\documentclass[12pt, letterpaper, oneside]{article}

\usepackage[T2A]{fontenc}
\usepackage[utf8]{inputenc}
\usepackage[english]{babel}
\usepackage{graphicx} % Required for inserting images
\usepackage{setspace}

\usepackage{xcolor}
\usepackage{amsmath}
\usepackage{amssymb}
\usepackage{bbm}
\usepackage[left=30mm, top=30mm, right=30mm, bottom=30mm]{geometry}
\newcommand{\RNumb}[1]{\uppercase\expandafter{\romannumeral #1\relax}}
\DeclareMathOperator*{\esssup}{ess\, sup}
\newtheorem{lemma}{Lemma}
\newtheorem{theorem}{Theorem}

\title{On polynomial recurrence property of  ``Markov-up'' processes}
\author{D.O. Kalikaeva\footnote{Moscow State University \& Institute for Information Transmission Problems, Moscow, Russia; email: diana.kalikaeva@math.msu.ru}}

\date{}

\begin{document}

\maketitle
\begin{abstract}
This work is a continuation of \cite{DK23}. The object of study is ``Markov-up processes'' on $\mathbb Z_+$ and the moment of downcrossing a certain barrier. The processes considered in this paper differ from Markov ones by the presence of a memory in certain parts of the trajectory. In our previous paper \cite{DK23} exponential recurrence conditions were established. In this paper polynomial recurrence properties are considered  under certain new assumptions.

\medskip


\medskip


\end{abstract}
\section{Introduction}
Let us recall that about two and a half decades ago, professor Alexander Dmitrievich Solovyev, in a personal conversation with his students, outlined the idea of a process that behaves like Markov during periods of growth, but as soon as the process makes a jump at the bottom, it has a memory and its further behavior depends on the entire trajectory of a continuous fall to the considered moment. For a long time his ideas were not implemented, and only recently in \cite{VV22} apparently the first model of this type was proposed.

In the real world, complex phenomena can often be observed that resemble the ``domino effect'', where one event leads to a chain of subsequent events. Similar phenomena are observed in economics, ecology and other fields. Special mathematical models are required to model and analyze such interconnected systems. One of these models may be the ``Markov-up'' process considered in this work.

Consider a discrete process $X_n, n > 0$ on $\mathbb Z_+$ or on $\mathbb Z_{0, \bar N}$, where $\bar N \in \mathbb N$. Further, without loss of generality, we assume that the process is defined on $\mathbb Z_+$. It is assumed that the following property holds for the process: with some functions $\phi(i, j)$ and $\psi(i, \dots, j)$, $ i, j \in \mathbb{Z_+}$, and for any $n \in \mathbb N$ 
$$\mathsf P(X_{n+1} = j|\mathcal F_n, X_n \geq X_{n-1}) = \phi(X_n, j),$$
$$\mathsf P(X_{n+1} = j|\mathcal F_n, X_n \leq X_{n-1}) = \psi(X_{\zeta_n}, \dots, X_n, j).$$
Note that the last function depends on a random number of variables.

Next consider the stopping moment $\tau$, which characterizes the first moment when the process crosses downwards a certain level $N\in\mathbb Z_+$. Positive recurrence for process was recently studied in \cite{VV22}. Under certain assumptions, it was found that $\mathsf E_x \tau < x + C_1$, where $C_1 < \infty$. Exponential recurrence was studied later in \cite{DK23}. Under reinforced assumptions, it has been proved that  $\mathsf E_x e^{\alpha\tau} \le C_2 e^{\alpha x} $ for some $\alpha >0$, where $C_2 < \infty$.

The goal of the paper is to find assumptions under which the existence of finite polynomial moments for $\tau$ is guaranteed.

\section{The model and the assumptions}

Let us consider the process $ X_n,$ $ n\geq{0}$ on $\mathbb{Z_+}$ (or on ${\mathbb{Z}_{0,\widetilde{N}}} = \{0, \dots, \widetilde{N}\}$ for some $\widetilde{N} \in \mathbb{N}$, $ \widetilde{N}<{\infty}$) and define random variables for $N \geq{0}$

\begin{equation}\label{zeta}
\zeta_n: = \inf(k\leq{n}:\Delta X_i = X_{i+1} - X_i<{0},  \forall i = k, \dots, n-1),
\end{equation}
\begin{equation}\label{xi}
\xi_n: = \sup(k\geq{n}:\Delta X_i = X_{i+1} - X_{i} \geq{0},  \forall i = n, \dots, k-1)\vee n,
\end{equation}
\begin{equation}\label{chi}
\chi_n: = \sup(k\geq{n}:\Delta X_i = X_{i+1} - X_{i}<{0}, \forall i = n, \dots, k-1)\vee n,
\end{equation}
\begin{equation}\label{tau}
\tau := \inf(t\geq{0}:X_t\leq{N}).
\end{equation}
Also, let
\begin{align}\label{tf}
   \hat{X}_{i,n} := X_i \mathbbm{1} (\min(\zeta_n, n)\leq{i}\leq{n}),\quad\widetilde{\mathcal{F}_n} := \sigma(\zeta_n; \hat{X}_{i,n}: 0\leq{i}\leq{n}).
\end{align}

Note that $\widetilde{\mathcal{F}_n}$ is not a filtration. And for all $n \in \mathbb{Z_+} $ $\widetilde{\mathcal{F}_n}\subset \mathcal{F}_n$, where $(\mathcal{F}_n)_{n \in \mathbb{Z_+}}$ is the natural filtration associated to the process $(X_n)_{n \in \mathbb{Z_+}}$.  

As usual, all inequalities with conditional probabilities and conditional expectations are understood a.s.

The following assumptions are made:

\textbf{A1. Random memory depth:}
\textit {For any $ n \geq{0} $}
        
\begin{equation}\label{ctF}
\mathsf P(X_{n+1} = j | \mathcal{F}_n) = \mathsf P(X_{n+1} = j | \widetilde{\mathcal{F}_n}) ,
\end{equation}
This condition distinguishes the proposed process from Markov chain. If $n$ is such that $X_{n-1} \le X_{n}$, then $\widetilde{\mathcal{F}}_{n} = \sigma (X_{n})$ and it may be informally said that at this moment the process has a  Markov property. Yet, if $X_{n-1} > X_{n}$, then $\widetilde{\mathcal{F}}_{n} = \sigma (\zeta_{n}, X_{\zeta _{n}}, \dots, X_{n} )$ and the conditional distribution of the next value $X_{n+1}$ depends on some nontrivial part of the past; informally it may be said that at this moment the behavior of the process is not markovian.
    
\textbf{A2. Irreducibility (local mixing):}
\textit {For any $x \leq{N} $ and for all two states $ y = x, y = x+1$}
\begin{center}
$\mathsf P(X_{n+1} = y | \widetilde{\mathcal{F}_n}, X_n = x) \geq{\rho} >{0}.$
\end{center}
Note that $2\rho \le 1$. The assumption A2 will guarantee the irreducibility of the process in the extended state space where the process becomes Markov.
        
\textbf{A3. Recurrence-1:}
\textit {There exists $N \geq{0}$ such that for the conditional probability of a jump down the following is performed}
\begin{equation}\label{q0}
\mathsf P(X_{n+1}<{X_n}|\widetilde{\mathcal{F}_n}, N<{X_n})\geq{\kappa_0}>{0},
\end{equation}
$$\mathsf P(X_{n+1}<{X_n}|\widetilde{\mathcal{F}_n}, N<{X_n}<{X_{n-1}})\geq{\kappa_1}>{0},$$ 
\textit{etc., for all $ n\geq{m} $}
\begin{equation}\label{q}
\mathsf P(X_{n+1}<{X_n}|\widetilde{\mathcal{F}_n}, N<{X_n}<{\dots}<{X_{n-m+1}})\geq{\kappa_{m-1}}>{0},\forall m\geq{1}.
\end{equation}
Note that $\{\kappa_n\}_{n = 0} ^{\infty}$ is a non-decreasing sequence.

Denote $q = 1 - \kappa_0 < 1$. Then for the conditional probability of a jump up we have
\begin{center}
$ \mathsf P(X_{n+1}\ge{X_n}|\widetilde{\mathcal{F}_n}, N<{X_n})\le 1 - {\kappa_0} = q < 1$.
\end{center}

The next two assumptions differ from initial ones proposed in article \cite{VV22}.

\textbf{A4. Recurrence-2:}
\textit {It is assumed that the following series  converges $\forall m < \infty$ }

        $$\sum\limits_{i\geq 1} i^m (1 - \kappa_i) < \infty.$$

\textit {It is assumed that the following infinite product converges}
\begin{align}\label{prod}
\bar\kappa_{\infty} := \prod\limits_{i=0}^{\infty}\kappa_i > 0.
\end{align}
This assumption implies that $\kappa_i \to 1$ as $i\to\infty$. Given that $\{\kappa_i\}_{i = 0} ^{\infty}$ is monotonic, its convergence to one  may be interpreted as follows: the longer the process is falling down, the higher is the probability that it will continue to fall.

Let $\bar q = 1 - \bar\kappa_{\infty} < 1$. This is the upper bound for probability that in one
go the process will not reach $[0, N]$.
        
\textbf{A5. Polynomial moment of the value of jump up is limited:}
\textit {For any $ m < \infty$}
\begin{align}\label{emoment} 
M_{\alpha}:=\esssup_{\omega \in \Omega}\sup_{n}\mathsf{E}((X_{n+1}-X_n)_+ ^ m|\widetilde{\mathcal{F}_n})<{\infty}.
\end{align}

\section{Polynomial moments} 

Before proceeding to the proof of the main results, we will prove a number of auxiliary lemmas. 
\begin{lemma}\label{lemma1}
Under the assumption (A3) for any $ x>{N}$ and $m < \infty$
\begin{equation}
\mathsf{E}_{x}(\xi_n-n)^m\leq{ M_2:={q(1 + q)\over{(1-q)^3}}}.		
\end{equation}
\end{lemma}
\textbf{Proof.} 
For all ${i \geq{n}}$ let
\begin{center}
$e_i = \mathbbm{1}(X_{i+1}\geq{X_i})$,  $\bar{e}_i = \mathbbm{1}(X_{i+1}\textless{X_i})$,  $l^i_n = \bar{e}_i\prod\limits_{k=n}^{i-1}e_k$,  $\Delta{X_i} = X_{i+1} - X_i$.
\end{center}
We have, for all ${i \geq{n}}$
\begin{center}
$\mathsf{E}_{x}(e_i|X_i\textgreater{N}) = \mathsf{P}_{x}(X_{i+1}\geq{X_i}|X_i\textgreater{N}) $ 
$$= \mathsf{E}_{x}(\mathsf{P}_{x}(X_{i+1}\geq{X_i}| \widetilde{\mathcal{F}_i}, X_i\textgreater{N})| X_i\textgreater{N}) \leq{1-\kappa_0} = q.$$
\end{center}
Then almost surely
\begin{center}
$(\xi_n-n)^m =  \sum\limits_{k=0}^{\infty}k^m\bar {e}_{n+k}\prod\limits_{i=n}^{n+k-1}e_i = \sum\limits_{k=1}^{\infty}k^m\bar {e}_{n+k}\prod\limits_{i=n}^{n+k-1}e_i = \sum\limits_{k=1}^{\infty}k^m l^{n+k}_n.$
\end{center}
Using this representation, get
\begin{center}
$\mathsf{E}_{x}(\xi_n-n)^m = \mathsf{E}_{x}\sum\limits_{k=0}^{\infty}k^m\bar{e}_{n+k}
\prod\limits_{i=n}^{n+k-1}e_i \le \sum\limits_{k=1}^{\infty}k^m
\mathsf{E}_{x}\prod\limits_{i=n}^{n+k-1}e_i \le \sum\limits_{k=1}^{\infty}k^m q^k =: \widetilde M_2 < \infty.$
 \end{center}
 The series converges $\forall m$ on the Cauchy's root test: $\sqrt[k]{a_k} = \sqrt[k]{k^m q^k} \to q, k\to\infty$, since $q < 1$ (by assumption (А5)).
 
Lemma 1 is proved. \hfill $\square$ 
\begin{lemma}\label{lemma2}
Under the assumptions (A3) and (A4) $\forall$ $ x\textgreater{N} $ and $m < \infty$
\begin{equation}
\mathsf{E}_{x}(\chi_n-n)^m\mathbbm{1}(\chi_n\textless{\tau})\leq\sum\limits_{i=1}^{\infty} i^m (1-\kappa_i)=:M_3\textless{\infty}.
\end{equation}
\end{lemma}
\textbf{Proof.}
As in the previous lemma, in the same notation we have
\begin{center}
$(\chi_n-n)^m\mathbbm{1}(\chi_n\textless{\tau}) \leq{}  \sum\limits_{k=1}^{\infty} k^m e_{n+k}\mathbbm{1}(n+k-1\textless{\tau})\prod\limits_{i=n}^{n+k-1}\bar{e}_i$.
\end{center}
So, 
\begin{center}
$\mathsf{E}_{x}(\chi_n-n)^m\mathsf{1}(\chi_n\textless{\tau}) \le
\mathsf{E}_{x} \sum\limits_{k = 1}^{\infty} k^m e_{n+k}
\mathsf{1}(n+k-1\textless{\tau})\prod\limits_{i=n}^{n+k-1}\bar{e}_i
\le \sum\limits_{k = 1}^{\infty}k^m\mathsf{E}_{x}
\mathbbm{1}(n+k-1\textless{\tau})\prod\limits_{i=n}^{n+k-1}\bar{e}_i
\mathsf{E}_{x}(e_{n + k} | \Delta{X_i}\textless{0}, n\leq{}i\leq{}n+k-1)
\le \sum\limits_{k = 1} ^ {\infty} k^m \mathsf{E}_{x}
\mathbbm{1}(n+k-1\textless{\tau})\prod\limits_{i=n}^{n+k-1}\bar{e}_i
(1 - \kappa_k) \le \sum\limits_{k = 1}^{\infty} k^m(1 - \kappa_k) 
=: M_3 \textless \infty
$
\end{center}
Lemma 2 is proved. \hfill $\square$ 

\begin{lemma}\label{lemma3}
Under the assumptions  (A3) and (А5) for all $ x\textgreater{N}$ and $m <\infty$
\begin{equation}
   %%%%%%%%%%%%%%%%%%
\sup_{n,x}\mathsf{E}_{x}((X_{\xi_n}-X_n)_+)^m \le M_4 \textless{\infty}.
\end{equation}
\end{lemma}
\textbf{Proof.}
With the same notation we have
\begin{center}
$\mathsf{E}_{x} ((X_{\xi_n} - X_n)_+)^m = \mathsf E_x 
\sum\limits_{i = n + 1}^\infty l_n^i (X_i - X_n)^m = 
\mathsf{E}_{x} \sum\limits_{i = n + 1}^\infty l_n^i (X_i - X_n)^m.$       
\end{center}
\begin{center}
$\mathsf E_x l_n^i (X_i - X_n)^m = \mathsf E_x 
(l_n^i (\sum\limits_{j = n} ^ {i-1} \Delta X_j)^m) = 
\mathsf E_x (l_n^i\sum\limits_{j = n} ^ {i-1} \Delta X_j)^m \leq \mathsf E_x l_n^i (i - n)^{m-1} \sum_{j = n}^{i-1}\Delta X_j^m.
$
\end{center}
Jensen's inequality was used in the end. %%%%%%%%%%%%%%% корявое
Further, 

\begin{center}
$\mathsf E_x  \sum\limits_{j = n} ^ {i-1}(\Delta X_j l_n^i)^m =
\sum\limits_{j = n} ^ {i-1} \mathsf E_x l_n^i (\Delta X_j)^m,
$
\end{center}
\begin{center}
$\mathsf E_x l_n^i (\Delta X_j)^m = \mathsf E_x (\Delta X_j)^m \bar e_i
\prod \limits_{k = n}^{i - 1} e_k \le \mathsf E_x (\Delta X_j)^m
\prod \limits_{k = n}^{i - 1} e_k = \mathsf E_x\mathsf 
E_{\mathcal F_{j+1}} (\Delta X_j)^m\prod \limits_{k = n}^{j} e_k
\prod \limits_{k' = j + 1} ^ {i - 1} e_{k'} $
$=\mathsf E_x (\Delta X_j)^m \prod \limits_{k = n}^{j} e_k
\mathsf E_{\mathcal F_{j+1}} \prod \limits_{k' = j + 1} ^ {i - 1} e_{k'}
\stackrel{A3}{\le} q^{i - j - 1} \mathsf E_x (\Delta X_j) ^m \prod
\limits_{k = n} ^ {j } e_k $
$= q ^ {i - j - 1} \mathsf E_x \mathsf E_{\mathcal F_j} (\Delta X_j) ^m \prod
\limits_{k = n} ^ {j - 1} e_k e_j = \mathsf E_x\prod\limits_{k = n}^{j - 1} e_k
\mathsf E_{\mathcal {\widetilde F}_j}(\Delta X_j) ^m e_j \stackrel{\widetilde {A5}}{\le} 
\widetilde M_1 q^{i - j - 1}\mathsf E_x\prod\limits_{k = n}^{j - 1} e_k$
$\stackrel{A3}{\le} \widetilde M_1 q^{i - j - 1} q^{j - n} = 
\widetilde {M}_1 q^{i - n - 1}.
$
\end{center}
Therefore, 
\begin{center}
$\mathsf{E}_{x} ((X_{\xi_n} - X_n)_+)^m = \sum\limits_{i = n + 1} ^\infty (i - n)^{m - 1} \sum\limits_{j = n} ^{i - 1} \widetilde {M}_1 q^{i - n - 1} = \sum\limits_{i = n + 1} ^\infty (i - n)^{m} \widetilde {M}_1 q^{i - n - 1}$
$=\widetilde {M}_1 q \sum\limits_{i = 1} ^\infty i^{m} q^{i} =: M_4 < \infty$
\end{center}
Lemma 3 is proved. \hfill $\square$

\medskip 
Now, we can proceed to the main theorem. Let us recall,
\begin{center}
$\tau := \inf(t\geq{0}:X_t\leq{N})$.
\end{center}  
\begin{theorem}\label{teorema}
Under the assumptions (A1), (A3)-(А5)  $\forall m$ there exist constants $C_1 < \infty$, $C_2 < \infty$ such that 
\begin{equation}
\mathsf{E}_{x} \tau ^ m  \leq C_1 (C_2 +x^m ).
\end{equation}
\end{theorem}
\textbf{Proof.}
If $x\leq{N} $, then $\tau = 0$, this case is trivial, therefore, we assume that $x\textgreater{N}$.
Consider the events:
\begin{center}
$ A_i = \{$exactly $i-1$ unsuccessful attempts to descend to the floor $[0, N],$ attempt no. $ i$ is successful $\}, i\geq{1},$\\
    
\medskip

$ B_j = \{$attempt no. $ j$ to reach $[0, N]$ is unsuccessful$\}, j\geq{1},$ \\           
\medskip

$ B_j^c = \{$attempt no. $ j$ to reach $[0, N]$ is successful$\}, j\geq{1}.$
\end{center}
Note that (according to the assumption А3) the probability of any unsuccessful attempt to cross the floor (that is, event $ B_j$) is less then $\bar q$. In this notations for event $A_i$ the following is valid: $\tau = T_i$, $A_i = (\cap_{j=1}^{i-1} B_j)\cap B_i^c$, $\mathsf P(A_i) \leq{\bar {q}^{i-1}}$.
\medskip

\textbf{Case \RNumb{1}}: at $t=0$ the process is falling down.
Let us define stopping times:
\begin{center}
$t_0 = T_0 = 0,$ $ T_1 = \chi_{t_0},$ $t_1 = \xi_{T_1},$ $ T_2 = \chi_{t_1},$ $t_2 = \xi_{T_2},$ $ T_3 = \chi_{t_2}, \dots$
\end{center}
So, from $t_{i-1}$ to $ T_i$ the process is continuously falling, at $T_i$ the fall is replaced by growth, and up to $ t_i$ the process continuously runs up. The process will change its behavior almost surely  finitely many  times  until it reaches the set $[0, N]$.
Note that $T_i-t_{i-1}\leq{X_{t_{i-1}}}$ and $B_j \in \mathcal F_{T_j}$.        
Let us estimate
\begin{center}
$\mathsf E_x  \tau^m = \sum\limits_{i\ge1}\mathsf E_x  \tau^m 
\mathbbm{1}(A_i) = \sum\limits_{i\ge1}\mathsf E_x T_i^m 
\mathbbm{1}(A_i) = \sum\limits_{i\ge1} \mathsf E_x T_i^m \mathbbm{1}
((\cap_{j=1}^{i-1}B_j)\cap B_i^c)$
$ = \sum\limits_{i\ge1}\mathsf E_x
(\prod\limits_{j=1}^{i-1}\mathbbm{1}(B_j)) \mathbbm{1}(B_i^c) T_i^m 
= \sum\limits_{i\ge1} \mathsf E_x\mathsf E_{\mathcal F_{t_{i-1}}}
(\prod\limits_{j=1}^{i-1}\mathbbm{1}(B_j)) \mathbbm{1}(B_i^c) T_i^m$
$= \sum\limits_{i\ge1} \mathsf E_x 
(\prod\limits_{j=1}^{i-1}\mathbbm{1}(B_j))
\mathsf E_{\mathcal F_{t_{i-1}}} \mathbbm{1}(B_i^c) T_i^m 
= \sum\limits_{i\ge1} \mathsf E_x (\prod\limits_{j=1}^{i-1}
\mathbbm{1}(B_j))\mathsf E_{\mathcal F_{t_{i-1}}} 
\mathbbm{1}(B_i^c) (T_i - t_{i-1} + t_{i-1})^m $
$\le 2^{m-1}\sum\limits_{i\ge1} \mathsf E_x (\prod\limits_{j=1}^{i-1}
\mathbbm{1}(B_j))\mathsf E_{\mathcal F_{t_{i-1}}} 
\mathbbm{1}(B_i^c) (T_i - t_{i-1})^m 
+ 2^{m-1}\sum\limits_{i\ge1} \mathsf E_x (\prod\limits_{j=1}^{i-1}
\mathbbm{1}(B_j))\mathsf E_{\mathcal F_{t_{i-1}}} 
\mathbbm{1}(B_i^c) t_{i-1}^m$
$\le 2^{m-1}\sum\limits_{i\ge1} \mathsf E_x (\prod\limits_{j=1}^{i-1}
\mathbbm{1}(B_j))\mathsf E_{\mathcal F_{t_{i-1}}} 
\mathbbm{1}(B_i^c) X_{t_{i-1}}^m 
+ 2^{m-1}\sum\limits_{i\ge1} \mathsf E_x (\prod\limits_{j=1}^{i-1}
\mathbbm{1}(B_j))\mathsf E_{\mathcal F_{t_{i-1}}} 
\mathbbm{1}(B_i^c) t_{i-1}^m$
$2^{m-1}\sum\limits_{i\ge1} \mathsf E_x (\prod\limits_{j=1}^{i-1}
\mathbbm{1}(B_j))\mathbbm{1}(B_i^c) X_{t_{i-1}}^m 
+ 2^{m-1}\sum\limits_{i\ge1} \mathsf E_x (\prod\limits_{j=1}^{i-1}
\mathbbm{1}(B_j))\mathbbm{1}(B_i^c) t_{i-1}^m$
\end{center}
We will evaluate separately $\mathsf E_x (\prod\limits_{j=1}^{i-1}
\mathbbm{1}(B_j))t_{i-1} ^m$ and $\mathsf E_x
(\prod\limits_{j=1}^{i-1}\mathbbm{1}(B_j))\mathbbm{1}(B_i^c) 
X_{t_{i-1}}^m.$

1) $$t_{i-1} = (t_{i-1} - T_{i-1}) + (T_{i-1} - t_{i-2})+\dots{}+(T_1-t_0)+(t_0-T_0).$$

$$ t_{i-1}^m = ((t_{i-1} - T_{i-1}) + (T_{i-1} - t_{i-2})+\dots{}
+(T_1-t_0)+(t_0-T_0))^m 
$$
$$\le 2^{m-1}[(\sum\limits_{j = 0}^{i-1}
(t_{j} - T_{j}))^m +(\sum\limits_{k = 0}^{i-2}(T_{k + 1} 
- t_k))^m]\le 2^{m-1}[i^{m-1}\sum\limits_{j = 0}^{i-1}
(t_j - T_j)^m + (i - 1)^{m-1}\sum\limits_{k = 0}^{i-2}(T_{k + 1} 
- t_k)^m].$$
We have
\begin{center}
$\mathsf E_x \prod\limits_{j=1}^{i-1}\mathbbm{1}(B_j)(t_{i-1} - 
T_{i-1})^m  = \mathsf E_x \mathsf E_{\mathcal F_{T_{i-1}}} 
\prod\limits_{j=1}^{i-1}\mathbbm{1}(B_j)(t_{i-1} - T_{i-1})^m $
$= \mathsf E_x\prod\limits_{j=1}^{i-1}\mathbbm{1}(B_j)
\mathsf E_{\mathcal F_{T_{i-1}}} (t_{i-1} - T_{i-1})^m 
\stackrel{lemma 1}{\le} M_2\mathsf E_x\prod\limits_{j = 1}^{i-1}
\mathbbm{1}(B_j)\le M_2 \bar q^{i-1}.$
\end{center}
Next,
\begin{center}
$\mathsf E_x \prod\limits_{j=1}^{i-1}\mathbbm{1}(B_j)
(T_{i - 1} - t_{i-2})^m = \mathsf E_x \mathsf E_{\mathcal F_{t_{i-2}}}
\prod\limits_{j=1}^{i-1}\mathbbm{1}(B_j)(T_{i - 1} - t_{i-2})^m$
$= \mathsf E_x \prod\limits_{j=1}^{i-2}\mathbbm{1}(B_j)
\mathsf E_{\mathcal F_{t_{i-2}}}\mathbbm{1}(B_{i-1})
(T_{i - 1} - t_{i-2})^m\stackrel{lemma2}{\le} M_3
\mathsf E_x \prod\limits_{j=1}^{i-2}\mathbbm{1}(B_j)\le
M_3 \bar q^{i-2}.$
\end{center}
Consider $k \in \mathbb N$ such that $1\le k\textless i $. Then we have, 
\begin{center}
$ \mathsf E_x \prod\limits_{j=1}^{i-1}\mathbbm{1}(B_j)(t_{k-1} - 
T_{k-1})^m  = \mathsf E_x \mathsf E_{\mathcal F_{t_{k-1}}} 
\prod\limits_{j=1}^{i-1}\mathbbm{1}(B_j)(t_{k-1} - T_{k-1})^m 
=  \mathsf E_x \prod\limits_{j=1}^{k-1}\mathbbm{1}(B_j)
(t_{k-1} - T_{k-1})^m\mathsf E_{\mathcal F_{t_{k-1}}} 
\prod\limits_{j'=k}^{i-1}\mathbbm{1}(B_{j'}) \le \bar q^{i-k}
\mathsf E_x \prod\limits_{j=1}^{k-1}\mathbbm{1}(B_j)
(t_{k-1} - T_{k-1})^m \le \bar q^{i-k} M_2\bar q^{k-1} = M_2\bar q^{i-1}.$
\end{center}
Now consider $l \in\mathsf N$ such that $2\le l \leq i$. Then 
\begin{center}
$\mathsf E_x \prod\limits_{j=1}^{i-1}\mathbbm{1}(B_j)
(T_{l - 1} - t_{l-2})^m = \mathsf E_x \mathsf E_{\mathcal F_{T_{l-1}}}
\prod\limits_{j=1}^{i-1}\mathbbm{1}(B_j)(T_{l - 1} - t_{l-2})^m$
$= \mathsf E_x \prod\limits_{j=1}^{l-1}\mathbbm{1}(B_j) \mathsf E_x 
(T_{l - 1} - t_{l-2})^m \mathsf E_{\mathcal F_{T_{l-1}}}
\prod\limits_{j=l}^{i-1}\mathbbm{1}(B_j)(T_{l - 1} - t_{l-2})^m$
$\le \bar q^{i-l}  \mathsf E_x \prod\limits_{j=1}^{l-1}
\mathbbm{1}(B_j) \mathsf E_x (T_{l - 1} - t_{l-2})^m \le\bar q^{i-l}
M_3 \bar q^{l-2} = \bar q^{i-2} M_3.$
\end{center}
\begin{center}
$\mathsf E_x (\prod\limits_{j=1}^{i-1}\mathbbm{1}(B_j))t_{i-1}^m
\le (2i)^{m-1} M_2 \sum\limits_{j = 0}^{i-1}\bar q^{i-1}
+ (2(i-1))^{m-1} M_3\sum\limits_{k = 0}^{i-2}\bar q^{i-2}$
$= 2^{m-1}i^m M_2\bar q^{i-1} + 2^{m-1}(i-1)^m M_3\bar q^{i-2}.$
\end{center}
2) Note that $X_{t_i} \geq X_{T_i}$ and $X_{T_0} = X_{t_0}=x$. Hence, 
\begin{center}
$X_{t_{i-1}} = (X_{t_{i-1}}- X_{t_{i-2}}) + (X_{t_{i-2}}- X_{t_{i-3}})+\dots+(X_{t_{1}}- X_{t_{0}} ) +(X_{t_0} - x) + x$
$\leq{}(X_{t_{i-1}}- X_{T_{i-1}}) + (X_{t_{i-2}}- X_{T_{i-2}})+\dots+(X_{t_{1}}- X_{T_{1}} ) +(X_{t_0} - X_{T_0}) + x,$
\end{center}
$X_{t_{i-1}}- X_{T_{i-1}} \ge 0$ for all $i$
\begin{center}
$X_{t_{i-1}}^m \le ((X_{t_{i-1}}- X_{T_{i-1}}) + (X_{t_{i-2}}- X_{T_{i-2}})+\dots+(X_{t_{1}}- X_{T_{1}} ) +(X_{t_0} - X_{T_0}) + x)^m$
$\le 2^{m-1}(\sum\limits_{j = 0}^{i-1}(X_{t_{j}} - X_{T_{j}}))^m + 2^{m-1}x^m
\le (2i)^{m-1} \sum\limits_{j = 0}^{i-1}(X_{t_{j}} - X_{T_{j}})^m + 2^{m-1}x^m.$
\end{center}
Then we estimate, 
\begin{center}
$\mathsf E_x (\prod\limits_{l=1}^{i-1}\mathbbm{1}(B_l))
\mathbbm{1}(B_i^c) \sum\limits_{j = 0}^{i-1}(X_{t_{j}} - X_{T_{j}})^m 
\le \sum\limits_{j = 0}^{i-1} \mathsf E_x (\prod\limits_{l=1}^{i-1}
\mathbbm{1}(B_l))(X_{t_{j}} - X_{T_{j}})^m.$
\end{center}
Consider $k \in \mathbb N$ such that $1\le k\le i - 1$. Then we have, 
\begin{center}
$\mathsf E_x (\prod\limits_{j=1}^{i-1}\mathbbm{1}(B_j))
(X_{t_k} - X_{T_k})^m = \mathsf E_x \mathsf E_{F_{t_k}}
(\prod\limits_{j=1}^{i-1}\mathbbm{1}(B_j))
(X_{t_k} - X_{T_k})^m$
$= \mathsf E_x (\prod\limits_{j=1}^{k}
\mathbbm{1}(B_j))(X_{t_k} - X_{T_k})^m\mathsf E_{F_{t_k}}
(\prod\limits_{j=k+1}^{i-1}\mathbbm{1}(B_j))\le\bar q^{i-k-1}
\mathsf E_x (\prod\limits_{j=1}^{k}\mathbbm{1}(B_j))
(X_{t_k} - X_{T_k})^m$
$=\bar q^{i-k-1}\mathsf E_x \mathsf E_{\mathcal F_{T_k}}
(\prod\limits_{j=1}^{k}\mathbbm{1}(B_j))(X_{t_k} - X_{T_k})^m
=\bar q^{i-k-1}\mathsf E_x (\prod\limits_{j=1}^{k}\mathbbm{1}(B_j))
\mathsf E_{\mathcal F_{T_k}}(X_{t_k} - X_{T_k})^m$
$\le \bar q^{i-k-1} M_4\mathsf E_x 
(\prod\limits_{j=1}^{k}\mathbbm{1}(B_j))\le \bar q^{i-k-1} M_4\bar q^k
=M_4\bar q^{i-1}.$
\end{center}

And for $k = 0$ we write, 

\begin{center}
$\mathsf E_x (\prod\limits_{j=1}^{i-1}\mathbbm{1}(B_j))
(X_{t_0} - X_{T_0})^m = \mathsf E_x \mathsf E_{F_{t_0}}
(\prod\limits_{j=1}^{i-1}\mathbbm{1}(B_j))
(X_{t_0} - X_{T_0})^m$
$= \mathsf E_x (X_{t_0} - X_{T_0})^m\mathsf E_{F_{t_0}}
(\prod\limits_{j=1}^{i-1}\mathbbm{1}(B_j))\le\bar q^{i-1}
\mathsf E_x 
(X_{t_0} - X_{T_0})^m$
$\le \bar q^{i-1} M_4.$
\end{center}

\begin{center}
$\mathsf E_x
(\prod\limits_{j=1}^{i-1}\mathbbm{1}(B_j))\mathbbm{1}(B_i^c) 
X_{t_{i-1}}^m \le \mathsf E_x 2^{m-1}x^m (\prod\limits_{j=1}^{i-1}\mathbbm{1}(B_j))
\mathbbm{1}(B_i^c) + (2i)^{m-1 }M_4\sum\limits_{j = 0}^{i-1}
\bar q^{i-1}$
$ = 2^{m-1}x^m \mathsf E_x (\prod\limits_{j=1}^{i-1}\mathbbm{1}(B_j))
\mathbbm{1}(B_i^c)+ 2^{m-1} i^m M_4\bar q^{i-1}.$
\end{center}

\medskip

\noindent
\textbf{Case \RNumb{2}}: at $t=0$ the process is going up.

Let us define stopping times:
\begin{center}
$T_0 = 0,$ $t_0 = \xi_{0},$ $ T_1 = \chi_{t_0},$ $t_1 = \xi_{T_1},$ $ T_2 = \chi_{t_1},$ $t_2 = \xi_{T_2},$ $ T_3 = \chi_{t_2}, \dots$
\end{center}

In this case we are using the same notation as in case Case \RNumb{1}. $ T_i$ is the end of the attempt no. $i$ to fall down, $ t_i$,  accordingly, is the moment when the growth is replaced by the fall. With probability one, the process will change its state only finitely many times  until it reaches the interval $[0, N]$.

In fact, the difference  between these cases is only in the meaning of the random value $t_0$. Since in Case \RNumb{1} $t_0 = T_0 = 0$, the terms $t_0 - T_0$ in step 1) and $X_{t_0} - X_{T_0}$ in step 2) are equal to zero and so, in fact, it was not necessary to take them into account in Case \RNumb{1}. Yet, the bounds obtained were a bit more general than required for the case, and, as a result, the estimates are valid in Case \RNumb{2}, too. 
%{\color{orange} Our goal is to obtain a universal estimate for both cases. Therefore, in Case \RNumb{1}, we did not make use of the fact that certain terms were zero, and we evaluated them too. As a result, all estimates obtained in Case \RNumb{1} are valid for Case \RNumb{2}. } {\color{blue} As we are considering these terms, all the estimates which were obtained in the previous part are valid for Case \RNumb{2}.} 
%{\color{red}(OK. A V ETOJ FRAZE CHTO SKAZANO? NE PONIMAYU, SKAZHITE YASNEE.)}

Now we may complete the proof of the theorem.

Since $\sum\limits_{i\ge1}\mathsf E_x (\prod\limits_{j=1}^{i-1}\mathbbm{1}(B_j))
\mathbbm{1}(B_i^c) = 1$, we estimate, 

\begin{center}
$\mathsf E_x  \tau^m \le 2^{m-1}\sum\limits_{i\ge1} \mathsf E_x (\prod\limits_{j=1}^{i-1}
\mathbbm{1}(B_j))\mathbbm{1}(B_i^c) t_{i-1} ^m $ 
$+ 2^{m-1}\sum\limits_{i\ge1} \mathsf E_x (\prod\limits_{j=1}^{i-1}
\mathbbm{1}(B_j))\mathbbm{1}(B_i^c) X_{t_{i-1}}^m $
$\le 2^{m-1}\sum\limits_{i\ge1} ( 2^{m-1}i^m M_2\bar q^{i-1} + 2^{m-1}(i-1)^m M_3\bar q^{i-2})$
$+ 2^{2m-2}x^m + 2^{m-1}\sum\limits_{i\ge1}2^{m-1} i^m M_4\bar q^{i-1}$
$= 2^{2m-2}(x^m + (M_2 + M_3 + M_4\bar q)\sum\limits_{i\ge1} i^m \bar q^{i-1} ) < C_1 (x^m + C_2),$
\end{center}

where $C_1 < \infty$, $C_2 < \infty$. 

Theorem \ref{teorema} is proved. \hfill $\square$

\section{Acknowledgments}

The study was funded by the Theoretical Physics and Mathematics Advancement Foundation “BASIS”.

\end{document}